\magnification=\magstep1
\input amstex
\documentstyle{amsppt}
\catcode`\@=11 \loadmathfont{rsfs}
\def\mycal{\mathfont@\rsfs}
\csname rsfs \endcsname \catcode`\@=\active

\vsize=7.5in

\topmatter
\title On the paving size of a subfactor  \\ 
$\text{\it in memory of Vaughan Jones and Mihai Pimsner}$
  \endtitle
 
\author Sorin Popa  \endauthor

\rightheadtext{Paving size}

\affil     {\it  University of California, Los Angeles} \endaffil

\address Math.Dept., UCLA, Los Angeles, CA 90095-1555, USA \endaddress
\email  popa\@math.ucla.edu\endemail

\thanks Supported in part by NSF Grant DMS-1955812  and the Takesaki Endowed Chair at UCLA \endthanks

\abstract Given an inclusion of II$_1$ factors $N\subset M$ with finite Jones index, $[M:N]<\infty$, 
we prove that for any $F\subset M$  finite and $\varepsilon >0$, 
there exists a partition of $1$ with $r\leq \lceil 16\varepsilon^{-2}\rceil$ $\cdot \lceil 4 [M:N]\varepsilon^{-2}\rceil$ projections $p_1, ..., p_r\in N$ 
such that $\|\sum_{i=1}^r p_ixp_i - E_{N'\cap M}(x)\|\leq \varepsilon \|x-E_{N'\cap M}(x)\|$, $\forall x\in F$ (where $\lceil \beta \rceil$ denotes the least integer $\geq \beta$). 
We consider a series of related invariants for $N\subset M$, generically called {\it paving size}. 
\endabstract

\endtopmatter

\document

\heading Introduction \endheading

A result in ([P97]) shows that an inclusion of separable II$_1$ factors $N\subset M$ has the so-called {\it relative Dixmier  property}, 
$\overline{\text{\rm co}} \{uxu^* \mid u\in \Cal U(N)\}\cap N'\cap M\neq \emptyset$ (where the closure is here in operator norm), 
for all $x\in M$, if and only if its Jones index is finite, 
$[M:N]<\infty$. 

Thus, if $[M:N]<\infty$ then given any $x\in M$ and any $\varepsilon>0$, there exist unitary elements $u_1, ..., u_n\in \Cal U(N)$ 
such that $\|\frac{1}{n}\sum_{i=1}^nu_ixu_i^* - E_{N'\cap M}(x)\|\leq \varepsilon$. Using this recursively, it follows that if $[M:N]<\infty$ then  
for any $F\subset M$ finite and any $\varepsilon>0$ there exist 
$v_1, ..., v_m \in \Cal U(N)$ such that $\|\frac{1}{m}\sum_{i=1}^mv_ixv_i^* - E_{N'\cap M}(x)\|\leq \varepsilon$, $\forall x\in F$. 

We attempt to identify in this paper the optimal number $n$ of unitaries necessary to ``$\varepsilon$-flatten'' this way an element $x$ 
(more generally a finite set $F$), exploring its dependence on $\varepsilon$ and on $[M:N]$. Our main result establishes an upper bound 
of magnitude $n\leq 64[M:N]\varepsilon^{-4}$, valid for any finite set $F\subset (M)_1$, arbitrarily large. 

The corresponding $n$ unitaries $u_1, ..., u_n\in N$ that we construct  are in fact powers $v^k, 0\leq k \leq n-1$, of a period $n$ unitary element $v\in \Cal U(N)$.  
Since an averaging by such $\{v^k\}_k$ satisfies $\frac{1}{n}\sum_{k=0}^{n-1}v^kxv^{-k} = \sum_{i=1}^n p_ixp_i$, where $p_i\in \Cal P(N)$ 
is a partition of $1$ with spectral projections of $v$, our result gives also an upper bound for the minimal size of a partition of 1 with projections 
$p_1, ..., p_n \in N$ with the property that $\|\sum_{i=1}^n p_ixp_i-E_{N'\cap M}(x)\|\leq \varepsilon$, $\forall x\in F$.  More precisely we get the following: 

\proclaim{Theorem} Let $N\subset M$ be an inclusion of $\text{\rm II}_1$ factors with finite Jones index, $[M:N]<\infty$. For any $F\subset M$ finite and any $\varepsilon >0$, 
there exists a partition of $1$ with $r \leq \lceil 16\varepsilon^{-2}\rceil$ $\cdot \lceil 4 [M:N]\varepsilon^{-2}\rceil$ projections $e_1, ..., e_r$ in $N$ such that  

\vskip.05in 

$ \quad \quad \quad \|\sum_{i=1}^r e_i x e_i - E_{N'\cap M}(x)\| \leq \varepsilon \|x-E_{N'\cap M}(x)\|, \forall x\in F.$  
\endproclaim

If  $x\in M$ has zero expectation onto $N'\cap M$, then  an expression of the form $\sum_i p_i xp_i$, with   
$p_i\in \Cal P(N)$ a partition of $1$ with projections in $N$ that diminishes to $\varepsilon$ the operator norm of $x$ 
is called an $\varepsilon$-{\it paving} of $x$ over $N$. Taking minimal size $n$ of partitions that can $\varepsilon$-pave a given $x\in M$ 
(or $F\subset M$), then the supremum of such $n$ over all $x\in (M)_1$ (or over all $F\subset (M)_1$ finite), gives numerical invariants for $N\subset M$ 
that we generically call 
{\it paving size} of $N\subset M$. The above result gives  the upper bound $64[M:N]\varepsilon^{-4}$ for all such invariants. Their exact calculation  
is an interesting problem. We comment on this and other related questions  in 
Section 2 of the paper (see the definitions, remarks and Corollary 2.5 in that section). This includes a discussion of the $L^2$-version of paving size invariants, 
in Remark 2.9.  

To prove the above result we first use 
(Theorem in [P92]) to obtain a partition of $1$ with $n\leq 16\varepsilon^{-2}$ projections 
$f_i=(f_{i,m})_m \in N^\omega$ (where $N^\omega$ is the ultrapower of $N$ with respect to some non-principal ultrafilter on $\Bbb N$) 
such that $\{f_i\}_i$ is free independent to the given finite set 
$F\subset M\ominus (N'\cap M)$. By (3.5 in [PV15]), this implies $\|\sum_{i=1}^n f_ixf_i\|\leq \varepsilon/2$, $\forall x\in F$, and so   
for any $\delta>0$, which one can take  arbitrarily small independently of any other constants involved ($\delta< \varepsilon^2/(4[M:N] |F|)^2$ will do), 
there is $m$ large enough such that $\|(\sum_i f_{i,m}xf_{i,m})(1-q_x)\|\leq \varepsilon/2 +\delta$, 
where $q_x\in \Cal P(M)$ are projections of trace $\leq \delta$, $\forall x\in F$. Due to the finiteness of  Jones' basic construction 
algebra $\langle M, e_N\rangle$ ([J82]), $E_N(q_x)$ have supports $s(E_N(q_x))$ of trace $\leq [M:N]\tau(q)\leq [M:N]\delta$, so they are all supported by a 
projection $p=\vee_{x\in F} s(E_N(q_x))$ of trace $\leq [M:N] |F|\delta$, that's still very small. This leaves room to flatten $p$ by a partition in $N$ 
with $\leq 4[M:N]\varepsilon^{-2}$ many projections, to make it $\leq  \varepsilon^2/4[M:N]$ in norm. Combining the two partitions, and using a key trick  
from (page 147 of [P98]), relying on 
the [PP83]-inequality $E_N(x)\geq [M:N]^{-1}x, \forall x\in M_+$, we deduce that this final partition, which has  
 $\leq (16\varepsilon^{-2})(4[M:N]\varepsilon^{-2})$ many projections,  paves all $x\in F$ to $\varepsilon/2 + \varepsilon/2=\varepsilon$. 

\vskip.05in

{\it Acknowledgement}. Like in the proof of the relative Dixmier property for finite index inclusions in (A.1 in [P96]; the Theorem and Corollary 4.1 in [P97]; Theorem 3.1 in [P98]), 
an important ingredient in the proof of its quantitative version above is played by the characterization 
of the Jones index $[M:N]$  that Mihai Pimsner and I have discovered in our paper ([PP83], INCREST preprint 52/1983): if $N\subset M$ is an inclusion of II$_1$ factors 
then $\lambda=[M:N]^{-1}$ satisfies $E_N(x)\geq \lambda x, \forall x\in M_+$, with $\lambda=[M:N]^{-1}$  the best constant for which such inequality holds true, i.e., 
$[M:N]$ $=(\sup \{c\geq 0 \mid E_N(x)\geq cx, \forall x\in M_+\})^{-1}$. We were led to this 
``probabilistic'' characterization of  $[M:N]$ while trying to elucidate some intriguing questions 
emanating from Vaughan Jones amazing paper {\it Index for subfactors} [J82], a preprint of which he sent us in the Summer of 1982. The present paper is in memory 
of the exciting exchanges of ideas, mathematical discussions and collaborations I had with Vaughan and with Mihai over the years. 
It is terribly sad to lose so dear friends. They will be greatly missed.

\heading 1. Proof of the Theorem  \endheading

For notations and terminology used hereafter we send the reader to ([P13], [AP17]), for 
basics in II$_1$ factors to ([AP17]), for subfactor theory to ([J82]). 

We first recall a Kesten-type norm estimate from ([PV15]): 

\proclaim{Lemma 1.1} Let $P$ be a $\text{\rm II}_1$ factor, $F=F^*\subset (P)_1$ a self-adjoint set of trace $0$ contractions and $n\geq 1$. 
Assume $v\in P$ is a unitary element with $v^n=1$, $\tau(v^k)=0$, $1\leq k < n$, such that $\{v\}''$ is free independent to $F\cup F^*$, 
i.e., $\tau(x_0\Pi_{i=1}^mv^{k_i}x_i)=0$ for all $m\geq 1$, $x_1, ..., x_{m-1} \in F$, $x_0, x_m \in F\cup \{1\}$, and $1\leq k_1, k_2, ..., k_m \leq n-1$. 
Then $\|\frac{1}{n}\sum_{k=1}^nv^kxv^{-k}\|\leq 2\sqrt{n-1}/n$, $\forall x\in F$. Equivalently, if $p_1, ..., p_n$ denote the minimal projections in $\{v\}''\simeq L(\Bbb Z/n\Bbb Z)$, then 
$\|\sum_{k=1}^n p_k x p_k\|\leq 2\sqrt{n-1}/n$, $\forall x\in F$.  
\endproclaim
\noindent
{\it Proof}. The freeness condition between the set $F$ and the algebra $\{v\}''$ implies that for any $x\in F$ the set $\{v^kxv^{-k}\mid 0\leq k \leq n-1\}$ is $L$-free 
in the sense of (Definition 3.1 in [PV15]). Thus, by (Corollary 3.5 in [PV15]), we have $\|\sum_{k=1}^n v^{k-1}xv^{-k+1}\|\leq 2\sqrt{n-1}$. 
The proof in [PV15] is based on (Proposition 3.4 in [PV15]), which shows that any $L$-free set of contractions $\{x_1, ..., x_n\}\subset N$ can be dilated 
to an $L$-free set of unitaries $\{U_1, ..., U_n\}$ in a larger II$_1$ factor $\tilde{N}\supset N$.  
Thus, one has  $\|\sum_{i=1}^n x_i \|\leq \|\sum_{i=1}^n U_i\|$. 
But the $L$-free condition  
for a set of unitaries $\{U_1, ..., U_n\}$ amounts to $\{U_1^*U_i\}_{i=2}^n$ being free independent Haar unitaries, for which one has $\|1+\sum_{i=2}^n U_1^*U_i\| 
=2\sqrt{n-1}$ by Kesten's Theorem ([K59]). When applied to the $L$-free set $x_k=v^{k-1}xv^{-k+1}$, $1\leq k \leq n$, this entails  
$$
\| \sum_{k=1}^n v^{k-1}xv^{-k+1}\| \leq \|\sum_{i=1}^n U_i\|
=\|1+\sum_{i=2}^n U_1^*U_i\|= 2\sqrt{n-1}. 
$$
\hfill $\square$ 

\proclaim{Lemma 1.2} Let $N\subset M$ be an inclusion of $\text{\rm II}_1$ factors with $\text{\rm dim}(N'\cap M)<\infty$, $F\subset (M)_1$ a finite 
set of elements with $0$ expectation onto $N'\cap M$ and $n\geq 1$.  
Given any $\delta>0$ there exists a partition of $1$ with projections $p_1, ...., p_n \in N$ and projections $q_i\in p_iMp_i$ of trace $\tau(q_i)\leq \delta$, $1\leq i \leq n$, 
such that $\|p_ix(p_i-q_i)\|\leq 2\sqrt{n-1}/n + \delta$, $\forall x\in F$, $1\leq i \leq n$. 
\endproclaim
\noindent
{\it Proof}. Let $\omega$ be a non-principal ultrafilter on $\Bbb N$. 
By (Theorem [P92]; see also Theorem 0.1 in [P13]), there exists $v\in \Cal U(N^\omega)$ such that $v^n=1$, $\tau(v^k)=0$, $1\leq k < n$, and such 
that the algebra $\{v\}''$ is free independent to $F \cup F^*$. If $f_1, ..., f_n\in \Cal P(N^\omega)$ are the minimal projection of $\{v\}''$,  then by Lemma 1.1 we have 
$\|f_ixf_i\|\leq 2\sqrt{n-1}/n$, $\forall x\in F$, $1\leq i \leq n$. 

Let $f_i=(f_{i,m})_m$ with $f_{i,m}\in \Cal P(N)$ and $\sum_i f_{i,m}=1$, $\forall m$. Since $F$ is finite, given any $\delta'>0$ there exists $m$ large enough such 
that the spectral projection $e_{x,i}$ of $(f_{i,m}xf_{i,m})^*(f_{i,m}xf_{i,m})$ corresponding to the interval $[4(n-1)/n^2 + \delta', \infty)$ has 
trace satisfying $\tau(e_{x,i})\leq \delta'$. Thus, if $\delta'$ is sufficiently small then the projection $q_i=\vee_{x\in F} e_{x,i} \in \Cal P(M)$, 
which has trace majorized by  $\sum_{x\in F}\tau(e_{x,i}) \leq \delta'|F|$, satisfies $\tau(q_i)\leq \delta$. 

It follows that if we let $p_i=f_{i,m}$ and $q_i=\vee_{x\in F} e_{x,i}$, 
then $\tau(q_i) =\sum_{x\in F} \tau(e_{x,i}) \leq \delta$, $q_i\leq p_i$ and for each $x\in F$, $1\leq i \leq n$ 
we have the norm estimate 
$$
\|p_ix(p_i-q_i)\|=\|(p_i-q_i)x^*p_ix(p_i-q_i)\|^{1/2} \leq \|(p_i-e_{x,i})x^*p_ix(p_i-e_{x,i})\|^{1/2}
$$
$$
=\|(f_{i,m}x^*f_{i,m}xf_{i,m})(f_{i,m}-e_{x,i})\|^{1/2}\leq2\sqrt{n-1}/n + \delta. 
$$
\hfill $\square$ 

\proclaim{Lemma 1.3}  Let $N$ be a $\text{\rm II}_1$ factor. For $b\geq 0$ in $N$, denote by $s(b)$ its support projection. 
Let $F\subset N$ be a finite set and $0< \varepsilon \leq 1/2$. Assume 
$2(\sum_{x\in F} \tau(s(|x|))) < \varepsilon(\max \{\|x\| \mid x\in F\})^{-1}$. Let $m$ denote the least integer greater than or equal to $\varepsilon^{-1} 
\max \{\|x\| \mid x\in F\}$. Then there exists a partition of $1$ with $m$ projections $q_1, ..., q_m \in N$ 
such that $\|\sum_{j=1}^m q_j x q_j\| \leq \varepsilon$. 
\endproclaim
\noindent
{\it Proof}. Let $e=\vee_{x\in F} (l(x) \vee r(x))$, where $l(x), r(x)$ denote the left and respectively right support projections of $x$. 

The condition $2(\sum_{x\in F} \tau(s(|x|))) \leq  \varepsilon(\max \{\|x\| \mid x\in F\})^{-1}$ together with the condition $m$ satisfies, 
imply that there exists a partition of $1$ with projections $e_1, ..., e_m\in N$ of trace $1/m$ such that $e\leq e_1$. Let $v\in \Cal U(N)$ be a unitary 
element satisfying $v^m=1$ and $v^{k-1} e_1 v^{-k+1}=e_k$, $1\leq k \leq m$. 

Let $q_1, ..., q_m$ denote the minimal projections of the abelian 
$m$-dimensional von Neumann algebra $\{v\}''$, with $v=\sum_{k=1}^{m} \alpha^{k-1} q_k$, where $\alpha=\exp(2\pi i/n)$. 
Since all $x\in F$ are supported on $e$ and $v^kev^{-k}$ are mutually disjoint, it follows that 
$\|\frac{1}{m}\sum_k v^kxv^{-k}\| \leq \|x\|/m$, $\forall x\in F$, which by the given conditions  gives 
$$\|\frac{1}{m}\sum_k v^kxv^{-k}\| \leq \varepsilon, \forall x\in F.$$
Since $\frac{1}{m}\sum_{k=0}^{m-1} v^kxv^{-k} = \sum_{k=1}^m q_k x q_k$, we are done.  
\hfill $\square$ 

\proclaim{Lemma 1.4}  Let $N\subset M$ be a an inclusion of $\text{\rm II}_1$ factors with finite Jones index, 
$[M:N]<\infty$. If $q\in M$ is a projection then $\tau(s(E_N(q)))\leq \tau(q)[M:N]$ 
\endproclaim
\noindent
{\it Proof}. Let $M\subset M_1:=\langle M, e_N \rangle$ be the basic construction for $N\subset M$, with $e_N\in M_1$ denoting as usual 
the corresponding Jones projection. Thus, $M_1=\text{\rm sp}  Me_NM$, $[N, e_N]=0$,  $e_Nxe_N=E_N(x)e_N$  
and $\tau(e_Nx)=\lambda \tau(x)$, $\forall x\in $, where $\lambda=[M:N]^{-1}$. 

If $q\in M$ is a projection, then one has $e_Nqe_N=E_N(q)e_N$. Thus, $s(e_Nqe_N)=s(E_N(q))e_N$ with its trace being equal to 
$\lambda \tau(s(E_N(q)))$. This implies that 
$$
\tau(q)\geq \tau(s(qe_Nq))  = \tau(s(e_Nqe_N))=\lambda \tau(s(E_N(q)), 
$$
and thus $\tau(s(E_N(q)) \leq \lambda^{-1}\tau(q)=[M:N]\tau(q)$. 
\hfill $\square$ 

\vskip.05in 
\noindent
{\it Proof of the Theorem}. Replacing $F$ by $\{x-E_{N'\cap M}(x)/\|x-E_{N'\cap M}(x)\| \mid x\in F\setminus N'\cap M\}$, we may assume $F\subset (M \ominus N'\cap M)_1$. By Lemma 1.2, for any given integer $n$ and any $\delta'>0$, there exists a partition of $1$ with projections 
$p_1, ..., p_n$ in $N$ of trace $1/n$ such that for each $1\leq i\leq n$ we have a projection $q_i\in p_iMp_i$ satisfying $\tau(q_i)\leq \delta'$ and 
$$
\|p_ix(p_i-q_i)\| \leq (4(n-1)/n^2 + \delta')^{1/2}, \forall x\in F.  \tag 1
$$ 

If we denote $b_{i,x}=q_ix^*p_ixq_i \in p_iMp_i$, $x\in F, 1\leq i\leq n$, then $b_{i,x}\in (p_iMp_i)_1$ are positive elements of support $\leq q_i$. 
It follows that $0\leq E_N(b_{i,x})\leq p_i$ and by Lemma 1.4, its support has trace $\tau(s(E_N(b_{i,x})))\leq [M:N]\tau(q_i)$. 

By Lemma 1.3, 
given any integer $m\leq \tau(p_i)/\tau(q_i)$, there exists a partition of $p_i$ with $m$ projections $q^i_1, ..., q^i_m \in \Cal P(p_iNp_i)$ of trace $\tau(p_i)/m$, 
such that 
$$
\|\sum_{j=1}^m q^i_j E_N(b_{i,x})q^i_j \|\leq 1/m.  \tag 2
$$
Since by (Theorem 2.1 in [PP83]) we have $b\leq [M:N] E_N(b)$ for any $b\in M_+$, it follows that 
$$
\|\sum_j q^i_j b_{i,x} q^i_j \| \leq  [M:N] \|\sum_{j=1}^m q^i_j E_N(b_{i,x})q^i_j \|\leq [M:N]/m. \tag 3 
$$
But since $\phi_i: p_iMp_i \rightarrow p_iMp_i$  defined by $\Phi_i(y)=\sum_j q^i_j y q^i_j$, $y\in p_iMp_i$, 
is unital completely positive, by Kadison's inequality we have $\phi_i(y^*)\phi_i(y) \leq \phi_i(y^*y)$, 
$\forall y\in p_iMp_i$. Applying this to $y=p_ixq_i$ and using $(3)$ it follows that for each $x\in F$ and $1\leq i \leq n$ we have
$$
\|\sum_j q^i_j (p_ixq_i)q^i_j \| \leq \|\sum_j q^i_j  (q_ix^*p_ixq_i)q^i_j\|^{1/2}   \tag 4 
$$
$$
=\|\sum_j q^i_j b_{i,x} q^i_j \|^{1/2} \leq ([M:N]/m)^{1/2}. 
$$
Also, since $\phi_i$ are contractive, by $(1)$ we have for each $i$ the estimate
$$
\|\sum_j q^i_j(p_ix(p_i-q_i))q^i_j\| \leq   \|p_ix(p_i-q_i)\| \tag 5
$$ 
$$\leq   (4(n-1)/n^2 + \delta')^{1/2}, \forall x\in F. 
$$ 

This implies that   
the partition of $1$ with $r=nm$ projections $\{e_k\}_{k=1}^r = \{q^i_j \mid 1\leq i \leq n, 1\leq j \leq m\}$, which refines $\{p_i\}_i$,  
satisfies for all $x\in F$ the inequalities
$$
\|\sum_k e_k x e_k \| \leq  \|\sum_{i,j} q^i_j(p_ix(p_i-q_i))q^i_j\| + \|\sum_{i,j} q^i_j (p_ixq_i)q^i_j \|\tag 6
$$
$$
\leq \max_i  \|p_ix(p_i-q_i)\| + \max_i \|\sum_{j} q^i_j (p_ixq_i)q^i_j \|
$$
$$
 \leq (4(n-1)/n^2 + \delta')^{1/2} + ([M:N]/m)^{1/2}  
$$

If we now take $\delta'< 4/n^2$ and the integers $n, m$ so that $m\geq 4[M:N]\varepsilon^{-2}$, $n\geq 16\varepsilon^{-2}$, then 
$(4(n-1)/n^2 + \delta')^{1/2} + ([M:N]/m)^{1/2}  \leq \varepsilon/2 + \varepsilon/2=\varepsilon$, ending the proof of the Theorem. 
\hfill $\square$

\heading 2. Further remarks  \endheading

\noindent
{\bf Definition 2.1}. If $N\subset M$ is an inclusion of II$_1$ factors with finite index, 
then for any $F\subset M$ non-empty and $\varepsilon >0$ we denote by 
$\text{\rm n}(N\subset M; F, \varepsilon)$ the infimum over all $n $ for which there exists a partition of $1$ with 
projections  $p_1, ..., p_n \in N$ 
such that $\|\sum_{i=1}^n p_i x p_i - E_{N'\cap M}(x)\|\leq \varepsilon \|x-E_{N'\cap M}(x)\|$, $\forall x\in F$, with the usual convention that this infimum is equal to $\infty$ 
if there exists no such finite partition. 
We call $\text{\rm n}(N\subset M; F, \varepsilon)\in \Bbb N\cup\{\infty\}$ 
the  $\varepsilon$-{\it paving size of $F$ in $N\subset M$}. 

\vskip.05in

\noindent
{\bf Definition 2.2}. For each $k=1, 2, ...$, we denote $\text{\rm n}_k(N\subset M; \varepsilon)\overset{def}\to{=} 
\sup \{\text{\rm n}(N\subset M; F, \varepsilon) \mid F\subset M_h, |F|\leq k \}$, where $M_h=\{x\in M \mid x=x^*\}$. We also denote $\text{\rm n}_{\infty}(N\subset M; \varepsilon) \overset{def}\to{=} 
\sup \{\text{\rm n}(N\subset M; F, \varepsilon) \mid \emptyset \neq F\subset M$ finite$\}$. 
These numbers are obviously isomorphism invariants 
for $N\subset M$ and we generically refer to them as 
{\it paving size} of $N\subset M$. 

Specifically, $\text{\rm n}(N\subset M; \varepsilon)=\text{\rm n}_1(N\subset M; \varepsilon)$ is called 
the $\varepsilon$-{\it paving size of} $N\subset M$ and for each $2\leq k \leq \infty$,  $\text{\rm n}_{k}(N\subset M; \varepsilon)$ 
is called $(\varepsilon,k)$-{\it paving size of} $N\subset M$. 

\vskip.05in

Note that these quantities are increasing in $k$, with $\sup_{k\geq 1} \text{\rm n}_{k}(N\subset M; \varepsilon)= \text{\rm n}_{\infty}(N\subset M; \varepsilon)$. So by the Theorem 
they are all bounded by an order of magnitude $64[M:N]\varepsilon^{-4}$. Also, if $N\subset P \subset M$ 
is an intermediate subfactor, then $\text{\rm n}_{k}(N\subset P; \varepsilon)\leq \text{\rm n}_{k}(N\subset M; \varepsilon)$, $\forall 1\leq k\leq \infty$. 

This terminology and notations are inspired by the similar ones used for MASAs (maximal abelian $^*$-subalgebras) in factors, $A \subset M$,  
in relation to the Kadison-Singer type problems (see e.g., [PV15]). Notably, the term ``paving'' was coined in relation with the Kadison-Singer problem and   
seems suitable for these quantities. 

Note that if $p_1, ..., p_n\in N$ is a partition of $1$ with projections and we denote $v=\sum_{k=1}^n \alpha^{k-1}p_k$, where 
$\alpha=\exp(2\pi i/n)$, then for any $x\in  M$ we have $\sum_{k=1}^n p_k x p_k=\frac{1}{n} \sum_{k=0}^{n-1} v^k x v^{-k}$. 
Thus, any ``paving'' of $x\in M$ with $n$-projections in a subfactor $N$ of $M$  (or in a MASA $A$ of $M$) can be viewed as a ``Dixmier averaging'' of $x$ 
by $n$-unitaries in $N$ (resp. $A$). 

\vskip.05in

\noindent
{\bf Definition 2.3}. In the same spirit as the pavings, for an inclusion of factors $N\subset M$, a finite set $\emptyset \neq F\subset M$  
and $\varepsilon >0$, we define the quantity $\text{\rm D}(N\subset M; F, \varepsilon)$ to be the infimum over all $n$ for which there exist $u_1, ..., u_n \in \Cal U(N)$ 
such  that $\|\frac{1}{n}\sum_{i=1}^n u_i x u_i^* - E_{N'\cap M}(x)\|\leq \varepsilon \|x-E_{N'\cap M}\|$, $\forall x\in F$. Then similarly to the above notations, 
we let $\text{\rm D}_\infty(N\subset M; \varepsilon) = \sup \{\text{\rm D}(N\subset M; F, \varepsilon) \mid \emptyset \neq F \subset M$ finite$\}$, 
$\text{\rm D}_k(N\subset M; \varepsilon) = \sup \{\text{\rm D}(N\subset M; F, \varepsilon) \mid \emptyset \neq F \subset M_h, |F|\leq k\}$, for $1\leq k <\infty$. 

\vskip.05in 

We clearly have $\text{\rm D}(N\subset M; F, \varepsilon)\leq \text{\rm n}(N\subset M; F, \varepsilon)$, for any finite $F\subset M$. Also,  
$\text{\rm D}_k(N\subset M; \varepsilon)\leq \text{\rm n}_k(N\subset M; \varepsilon)$,  for any $1\leq k \leq \infty$. So the Theorem implies that  
for any subfactor of finite index $N\subset M$, these quantities are all finite, in fact bounded by the order of magnitude $64[M:N]\varepsilon^{-4}$. 
Like the $\text{\rm n}_*(N\subset M;\varepsilon)$-quantities, 
they are all isomorphism invariants for $N\subset M$. We'll still view them as {\it paving}-invariants for $N\subset M$, but with respect 
to averaging by unitaries, rather than by projections summing up to $1$.  Alternatively, we view  them as {\it optimal Dixmier averaging numbers for} $N\subset M$. 

In particular, for a single II$_1$ factor $N$ and $1\leq k \leq \infty$ 
we have $\text{\rm n}_k(N;  \varepsilon)\overset{def}\to{=}\text{\rm n}_k(N\subset N; \varepsilon) \leq 64\varepsilon^{-4}$.  
Consequently,   $\text{\rm D}_k(N;  \varepsilon)\overset{def}\to{=}\text{\rm D}_k(N\subset N; \varepsilon) \leq 64\varepsilon^{-4}$ as well. 

Dixmier's classical 
{\it averaging theorem} (see Ch. III, Sec. $5$ in [D57]) amounts to $\text{\rm D}(N;\varepsilon):=\text{\rm D}_1(N;\varepsilon)< \infty$. His proof  actually shows that 
$\text{\rm D}(N;\varepsilon)\leq \lceil \varepsilon^{-c} \rceil$, where  $c=\log_{3/2}2=\frac{\ln 2}{\ln 3 - \ln 2}\approx 1.7095 < 2$. If $F$ is a finite set 
of $k$ selfadjoint elements, then by applying consecutively Dixmier's theorem $k$ many times, one obtains the estimate $\text{\rm D}_k(N;\varepsilon) \leq 
\lceil \varepsilon^{-c} \rceil^k$, which thus depends on $k$ and gives no bound for  $\text{\rm D}_\infty(N;\varepsilon)$. 
So Dixmier's proof gives better upper bounds for $\text{\rm D}_k(N; \varepsilon)$ if $k=1, 2$, but  a (exponentially)  worse bound for $k\geq 3$,  
with no bound for $k=\infty$. 

It would be interesting to improve the upper bound for the paving size $\text{\rm n}_k(N\subset M; \varepsilon)$, especially for $k=1, k=\infty$,  
as well as for the constants $\text{\rm D}_k(N\subset M; \varepsilon)$. In particular, to determine if the order of magnitude $\varepsilon^{-4}$ is optimal or can 
be lowered. Equally interesting would be to obtain some sharp lower bounds. Ideally, one would like to have exact calculation of $\text{\rm n}_*(N\subset M; \varepsilon)$  
or  $\text{\rm D}_*(N\subset M;\varepsilon)$, for some concrete subfactors $N\subset M$ of finite index. This seems quite challenging even for $N=M$! 

Another interesting problem is  
to determine whether these invariants only depend on the index $[M:N]$ (respectively, only on the standard invariant  $\Cal G_{N\subset M}$).  

One can provide a (rather weak!) estimate for the lower bound of the paving size constants from the following simple observation for single II$_1$ factors:  

\proclaim{Lemma 2.4} Let $N$ be a $\text{\rm II}_1$ factor. Let $x\in N_+$  be so that $\|x\|=1$. If $u_1, ..., u_n\in \Cal U(N)$ are so 
that $\|\frac{1}{n}\sum_i u_ixu_i^*-\tau (x)1\|\leq \varepsilon$, then $n\geq (\tau(x)+\varepsilon)^{-1}$. 
In particular, if $x=q\in \Cal P(N)$ is a non-zero projection, 
then $\text{\rm D}(N; q, \varepsilon)\geq (\tau(q)+\varepsilon)^{-1}$. 
\endproclaim 
\noindent
{\it Proof}. Since $\|\frac{1}{n}\sum_i u_ixu_i^*-\tau (x)1\|\leq \varepsilon$ and $x\geq 0$, we have 
$ (\tau (x) + \varepsilon)1 \geq \frac{1}{n}\sum_i u_ixu_i^* \geq \frac{1}{n} u_1xu_1^*$, 
so by taking norms we get $(\tau (x) + \varepsilon) \geq \| \frac{1}{n}\sum_i u_ixu_i^*\|\geq \frac{1}{n} \|x\| = \frac{1}{n},$
implying that $n \geq (\tau(x)+\varepsilon)^{-1}$. 
\hfill $\square$ 

\vskip.05in 

Taking $\tau(q) \rightarrow 0$ in Lemma 2.4 we get the lower bound $\varepsilon^{-1}$ for the paving size of a single II$_1$ factor, 
and hence for any 
inclusion of II$_1$ factors. Combining with the Theorem and the above remarks, we thus  get: 

\proclaim{Corollary 2.5} If $N\subset M$ is an inclusion of $\text{\rm II}_1$ factors with finite Jones index then, with  
the above notations, we have for any $\varepsilon >0$ the estimates 
\vskip.05in 
\noindent
$ \ \varepsilon^{-1} \leq \text{\rm D}(N\subset M; \varepsilon) 
\leq \text{\rm n}(N\subset M; \varepsilon) \leq \text{\rm n}_\infty(N\subset M; \varepsilon)\leq \lceil 16\varepsilon^{-2}\rceil \cdot \lceil 4 [M:N]\varepsilon^{-2}\rceil $
\endproclaim

The invariants $\text{\rm n}_*(N\subset M; \varepsilon)$, $\text{\rm D}_*(N\subset M; \varepsilon)$ can also be viewed as measuring 
how efficient one can ``flatten'' the elements in $M_+$ by averaging/paving with unitaries (or partitions with projections) in $N$. 
Two other quantities that measure such phenomena are the following: 

\vskip.05in 
\noindent
{\bf Definition 2.6}. Let $N\subset M$ be an inclusion of II$_1$ factors with finite index. Recall from   
(Corollary 3.1.9 in [J82])  that there exist projections $e\in M$ satisfying $E_N(e)=[M:N]^{-1}1$  
and that by (Corollary 1.8 in [PP83]) any two such projections are conjugate by a unitary in $N$. Thus, the quantity  
$\text{\rm d}(N\subset M; \varepsilon)\overset{def}\to{=} \text{\rm n}(N\subset M; e, \varepsilon)$, 
where $e$ is such a ``Jones projection'', is well defined and it is obviously an isomorphism invariant for $N\subset M$. 
One clearly has $\text{\rm d}(N\subset M; \varepsilon) \leq \text{\rm n}(N\subset M; \varepsilon)$. 

In a related vein, we define the invariant $\text{\rm d}_{\text{\rm ob}}(N\subset M)$  for a subfactor of finite index $N\subset M$  as the 
infimum of $\|\sum_j m_j^*m_j\|$ over all orthonormal basis $\{m_j\}_j$ of $N \subset M$  (as defined in Section 1 of [PP83]). 

\vskip.05in

Since for any orthonormal basis $\{m_j\}_j$ one has $\lambda \sum_j m_j m_j^* = 1$, where $\lambda=[M:N]^{-1}$,   
 (cf. Proposition 1.3 in [PP83]), it follows that $1=\tau(\lambda \sum_jm_jm_j^*)=\lambda \tau(\sum_j m_j^*m_j)$, hence $\|\sum_j m_j^*m_j\| \geq 
 \lambda^{-1} =[M:N]$. Thus, one has $\text{\rm d}_{\text{\rm ob}}(N\subset M)\geq [M:N]$. On the other hand, one can take the orthonormal basis 
$\{m_j\}_j$ so that $m_1=1$ and  so that for all but possibly one $m_j$ to have $E_N(m_j^*m_j)=1$, which by (Proposition 2.1 in [PP83]) implies $m_j^*m_j \leq [M:N] 1$. Thus 
$\|\sum_j m_j^*m_j\| \leq \sum_j \|m_j^*m_j\| \leq 1 + [M:N] (\lceil [M:N]\rceil -1).$

We have thus proved the following 

\proclaim{Proposition 2.7} If $N\subset M$ is an inclusion of $\text{\rm II}_1$ factors with finite Jones index then, with  
the above notations, we have the estimates 

\vskip.05in 

$\quad  \quad \quad \quad [M:N] \leq \text{\rm d}_{\text{\rm ob}}(N\subset M) \leq 1 + [M:N] (\lceil [M:N]\rceil -1).$
 \endproclaim

\vskip.05in 
\noindent
{\bf Remark 2.8}. The paving size invariants can be defined for an arbitrary inclusion of factors (not necessarily II$_1$), $\Cal N \subset \Cal M$, with exactly same formal definitions. 
If one has an expectation $\Cal E:\Cal M\rightarrow \Cal N$ with finite Pimsner-Popa index, i.e., if $\Cal E(x)\geq \lambda x$, $\forall x\in \Cal M_+$, 
for some $\lambda>0$, and one denotes by $\text{\rm Ind}(\Cal E)$ the inverse $\lambda^{-1}$ of the best constant $\lambda$ satisfying the inequality, 
then the main result in ([P97]) shows that Ind$(\Cal E)<\infty$ implies $\text{\rm n}(\Cal N\subset^{\Cal E} \Cal M; F, \varepsilon)<\infty$, for any finite set $F\subset \Cal M$. 
We leave it to the interested reader to adapt the proof of the Theorem in this paper, combined with the proof of the relative Dixmier property for 
inclusions of properly infinite factors $\Cal N \subset^{\Cal E} \Cal M$ with $\text{\rm Ind}(\Cal E)<\infty$ in (Section 3 of [P97]), to get estimates 
for $\text{\rm n}_*(\Cal E; \varepsilon)$,  $\text{\rm D}_*(\Cal E; \varepsilon)$. 

\vskip.05in
\noindent
{\bf Remark 2.9}. One can consider exactly the same type of definitions as we did 
for $\text{\rm n}_*(N\subset M; \varepsilon)$, $\text{\rm D}_*(N\subset M; \varepsilon)$, where we replace the operator norm by the Hilbert norm-$\| \ \|_2$ given by the trace.   
We denote these invariants of a subfactor $N\subset M$ 
by $\text{\rm n}_k^{(2)}(N\subset M; \varepsilon)$, $\text{\rm D}_k^{(2)}(N\subset M; \varepsilon)$, $\text{\rm d}^{(2)}(N\subset M; \varepsilon)$, respectively,  
and refer to them generically as $L^2$-{\it paving size} of $N\subset M$ (inspired by terminology used in Section 3 of [P13]). 
These invariants may be easier to calculate, but less relevant of the properties of the inclusion $N\subset M$. Recall in this respect that for any 
inclusion of II$_1$ factors $N\subset M$, the subfactor $N$ contains a MASA 
$A\subset N$ such that $A'\cap M=A \vee (N'\cap M)$ (see e.g., Corollary 1.2.3 in [P16]), which by (Theorem 3.6 in [P13]) contains 
approximate 2-independent partitions of any size. Thus, for any $F\subset M\ominus (A\vee N'\cap M)$ finite, any $\delta>0$ and any 
$n\geq 1$, one can find a partition of $1$ with projections of trace $1/n$ in $A$, $p_1, ..., p_n\in \Cal P(A)$, such that $\|\sum_{i=1}^n p_ixp_i\|_2 
\approx_\delta n^{-1/2} \|x\|_2$, $\forall x\in F$. Thus, one has  the estimates 
$\text{\rm D}^{(2)}_k(N\subset M; \varepsilon) \leq \text{\rm n}^{(2)}_k(N\subset M; \varepsilon)  \leq 
\text{\rm n}^{(2)}_{\infty}(N\subset M; \varepsilon) \leq [\varepsilon^{-2}]+1$, for any $1\leq k \leq \infty$, for any $N\subset M$, without even assuming $[M:N]<\infty$.   

\vskip.05in
\noindent
{\bf Remark 2.10}. The most interesting  case of inclusions of factors $\Cal N\subset \Cal M$ is when they are {\it ergodic}, i.e., $\Cal N'\cap \Cal M=\Bbb C$. They 
correspond to the action $\Cal U(\Cal N)\curvearrowright^{\text{\rm Ad}} \Cal M$ being ergodic. A strengthening of ergodicity, 
called  {\it MV-ergodicity} ([P19]), 
requires that the $wo$-closure of the convex hull of $\{uxu^*\mid u\in \Cal U(\Cal N)\}$ intersects $\Cal N'\cap \Cal M=\Bbb C1$ 
(see also [P98] where this is called {\it weak relative Dixmier property}). Since 
$wo$ and $so$-closures coincide on bounded convex sets, it is equivalent to $\overline{\text{\rm co}}^{so}\{uxu^*\mid u\in \Cal U(\Cal N)\}\cap \Bbb C1\neq \emptyset$.  
For an inclusion of II$_1$ factors $N\subset M$ this amounts to a von Neumann type $L^2$-mean value ergodicity: 
$\forall x\in M$, $\forall \varepsilon >0$, $\exists u_1, ..., u_n\in \Cal U(N)$ 
such that $\|\frac{1}{n}\sum_{i=1}^n u_ixu_i^*-\tau(x)1\|_2\leq \varepsilon$.  Viewed from this perspective, 
Dixmier's averaging theorem states that for any single factor $\Cal N$, the action 
$\Cal U(\Cal N)\curvearrowright^{\text{\rm Ad}} \Cal N$ is $L^\infty$-MV ergodic, while the result in (A.1 in [P96], [P97]) shows that 
$\Cal U(N)\curvearrowright^{\text{\rm Ad}} M$ is $L^\infty$-MV ergodic for any ergodic inclusion  of II$_1$ factors $N\subset M$ with finite Jones index $[M:N]<\infty$  
(with the converse holding true when $N, M$ are separable, by Corollary 4.1 in [P97]). Our results in this paper can be viewed as  quantitative estimates  
of $L^\infty$-MV ergodicity for finite index inclusions.

\head  References \endhead

\item{[AP17]} C. Anantharaman, S. Popa: ``An introduction to II$_1$ factors'', \newline www.math.ucla.edu/$\sim$popa/Books/IIun-v13.pdf 

\item{[D57]} J. Dixmier: ``Les alg\'ebres d'operateurs sur l'espace Hilbertien (Alg\'ebres de von Neumann)'', Gauthier-Villars, Paris, 1957.

\item{[J82]} V. F. R. Jones: {\it Index for subfactors}, Invent. Math., {\bf 72}  (1983), 1-25. 

\item{[KS59]} R.V. Kadison, I.M. Singer: {\it Extensions of pure
states}, Amer. J. Math. {\bf 81} (1959), 383-400.

\item{[K59]} H. Kesten: {\it Symmetric random walks on groups}, Transactions of the AMS, {\bf 92} (1959), 336-354.

\item{[PP83]} M. Pimsner, S. Popa: {\it Entropy and index for subfactors}, Ann. Sci. Ecole Norm. Sup., {\bf 19} (1986), 57-106 (INCREST preprint No. 52/1983). 

\item{[P92]} S. Popa: {\it Free independent sequences in type $\text{\rm II}_1$ factors and related problems}, Asterisque, {\bf 232} (1995), 187-202.

\item{[P96]} S. Popa:  {\it Some properties of the symmetric enveloping algebras
with applications to amenability and property T},
Documenta Mathematica, {\bf 4} (1999), 665-744.

\item{[P97]} S. Popa: {\it The relative
Dixmier property for inclusions of von Neumann
algebras of finite index}, Ann. Sci. Ec. Norm. Sup.,
{\bf 32} (1999), 743-767. 

\item{[P98]} S. Popa: {\it On the relative Dixmier property for
inclusions of C$^*$-algebras}, Journal of Functional Analysis, {\bf
171} (2000), 139-154.

\item{[P13]} S. Popa: {\it Independence properties in subalgebras of ultraproduct} II$_1$ {\it factors}, Journal of Functional Analysis 
{\bf 266} (2014), 5818Ð5846 (math.OA/1308.3982). 

\item{[P16]}  S. Popa:  {\it Constructing MASAs with prescribed properties}, Kyoto J. of Math,  {\bf 59} (2019), 367-397 (math.OA/1610.08945).

\item{[P19]} S. Popa: {\it On ergodic embeddings of factors}, Communications in Mathematical Physics,  {\bf 384} (2021), 971-996 (math.OA/1910.06923).

\item{[PV15]}  S. Popa, S. Vaes: {\it On  the optimal paving over MASAs in von Neumann algebras},  Contemporary Mathematics
Volume {\bf 671}, R. Doran and E. Park editors, American Math Society, 2016, pp 199-208 (math.OA/1507.01072).

\enddocument